\documentclass[12pt]{article}
\usepackage{amssymb}
\begin{document}

\title{\bf $G$-subsets and $G$-orbits of $Q^*(\sqrt{n})$
under action of the Modular Group.}

\author{M. Aslam Malik \thanks{malikpu@yahoo.com} and M. Riaz \thanks{mriaz@math.pu.edu.pk}\\
Department of Mathematics, University of the Punjab,\\
Quaid-e-Azam Campus, Lahore, Pakistan.}

\date{}
\maketitle
\begin{abstract}

It is well known that $G=\langle x,y:x^2=y^3=1\rangle$ represents
the modular group $PSL(2,Z)$, where $x:z\rightarrow\frac{-1}{z},
 y:z\rightarrow\frac{z-1}{z}$ are linear fractional
transformations. Let $n=k^2m$, where $k$ is any non zero integer
and $m$ is square free positive integer. Then the set
$$Q^*(\sqrt{n}):=\{\frac{a+\sqrt{n}}{c}:a,c,b=\frac{a^2-n}{c}\in
Z~\textmd{and}~(a,b,c)=1\}$$ is a $G$-subset of the real quadratic
field  $Q(\sqrt{m})$  \cite{R9}.  We denote
$\alpha=\frac{a+\sqrt{n}}{c}$ in $ Q^*(\sqrt{n})$ by
$\alpha(a,b,c)$. For a fixed integer $s>1$, we say that two
elements $\alpha(a,b,c)$, $\alpha'(a',b',c')$ of $Q^*(\sqrt{n})$
are $s$-equivalent if and only if $a\equiv a'(mod~s)$, $b\equiv
b'(mod~s)$ and $c\equiv c'(mod~s)$. The class $[a,b,c](mod~s)$
contains all $s$-equivalent elements of $Q^*(\sqrt{n})$ and
$E^n_s$ denotes the set consisting of all such classes of the form
$[a,b,c](mod~s)$.

In this paper we investigate proper $G$-subsets and $G$-orbits of
the set $Q^*(\sqrt{n})$ under the action of Modular Group $G$.\\

\end{abstract}



\section{Introduction}
An integer $m>0$ is said to be square free if its prime
decomposition contains no repeated factors. It is well known that
every irrational member of $Q(\sqrt{m})$ can be uniquely expressed
as $\frac{a+\sqrt{n}}{c}$, where $n=k^2m$ for some
integer $k$ and $a,~\frac{a^2-n}{c}$ and $c$ are relatively prime integers.\\

The set
$Q^*(\sqrt{n})=\{\frac{a+\sqrt{n}}{c}:a,c,b=\frac{a^2-n}{c}\in
Z~\textmd{and}~(a,b,c)=1\}$ is a proper $G$-subset of
$Q(\sqrt{m})$ \cite{R9}. If $\alpha=\frac{a+\sqrt{n}}{c}$ and
$\bar{\alpha}={\frac{a-\sqrt{n}}{c}}$ have different signs, then
$\alpha$ is called an ambiguous number. These ambiguous numbers
play an important role in the study of action of $G$ on
$Q(\sqrt{m})\cup\{\infty\}$, as $\textmd{Stab}_{\alpha}(G)$ are
the only non-trivial stabilizers and in the orbit $\alpha^G$,
there is
only one (up to isomorphism).\\
G. Higman (1978) introduced the concept of the coset diagrams for
the modular group $PSL(2,Z)$ and Q. Mushtaq (1983) laid its
foundation. By using the coset diagrams for the orbit of the
modular group $G=\langle x,y:x^2=y^3=1\rangle$ acting on  the real
quadratic fields
 Mushtaq \cite{R9} showed that for a fixed non-square positive
 integer $n$, there are only a finite number of ambiguous numbers
 in $Q^*(\sqrt{n})$, and that the ambiguous numbers in the coset
 diagram for the orbit $\alpha^G$ form a closed path and it is
 the only closed path contained in it. Let $C'=C\cup\{\pm\infty\}$ be
 the extended complex plane. The action of the modular group $PSL(2,Z)$
 on an imaginary quadratic field, subsets of $C'$, has been discussed in \cite{R2}.
   The action of the modular group on the real quadratic fields, subsets of $C'$,
     has been discussed in detail
 in \cite{R8}, \cite{R9} and \cite{R10}. The exact number of ambiguous
 numbers in $Q^*(\sqrt{n})$ has been determined in \cite{R11}, \cite{R6}
 as a function of $n$. The ambiguous length of an orbit $\alpha^G$ is the
 number of ambiguous numbers in the same orbit \cite{R11}, \cite{R6}.
The Number of Subgroups of $PSL(2,Z)$ when acting on
$F_p\cup\{\infty\} $ has been discussed in \cite{R14} and the
subgroups of the classical modular group
 has been discussed in \cite{R15}.  \\
 A classification of the elements $(a+\sqrt{p})/c,~b=(a^2-p)/c,$
 of $Q^*(\sqrt{p}),~p$ an odd prime, with respect to odd-even nature
 of $a,b,c$ has been given in \cite{R3}. M. Aslam Malik et al. \cite{R7}
 proved, by using the notion of congruence, that for each non-square
 positive integer $n>2$, the action of the group $G$ on a subset
 $Q^*(\sqrt{n})$ of the real Quadratic field $Q(\sqrt{m})$ is
 intransitive.\\
 If $p$ is an odd prime,then  $t~{\not\equiv}~0(mod~p)$ is said to
be a quadratic residue of $p$
if there exists an integer $u$ such that $u^2\equiv t(mod~p)$.\\
The quadratic residues of $p$ form a subgroup $Q$ of the group of
nonzero integers modulo $p$ under multiplication and
$|Q|=(p-1)/2$.  \cite{R1}\\
\textbf{Lemma 1.1}  If $v_1,v_2\in Q$, $n_1, n_2{\not\in}~Q$
($v_1,v_2$ are quadratic residues, and $n_1,n_2$ are
quadratic non-residue, Then\\
(a) $n_1v_1$ is a quadratic non-residue.\\
(b)$n_1n_2$ is a quadratic residue.\\
(c)$v_1v_2$ is a quadratic residue.\\
 In the sequel, $q.r$ and $q.nr$ will stand for quadratic
residue and quadratic non-residue respectively. The Legendre
symbol $(a/p)$ is defined as $1$ if $a$ is a quadratic residue of
$p$ otherwise it is
defined by $-1$.  \cite{R1}\\
We denote the element $\alpha=\frac{a+\sqrt{n}}{c}$of
$Q^*(\sqrt{n})$ by $\alpha(a,b,c)$ and say that two elements
$\alpha(a,b,c)$ and $\alpha'(a',b',c')$ of $Q^*(\sqrt{n})$ are
$s$-equivalent (and write $\alpha(a,b,c)\thicksim_s
\alpha'(a',b',c')$ or $\alpha\sim_s \alpha')$ if and only if
$a\equiv a'(mod~s)$, $b\equiv b'(mod~s)$ and $c\equiv c'(mod~s)$.
Clearly the relation $\thicksim_s$ is an equivalence relation, so
for each integer $s>1$, we get different equivalence classes
$[a,b,c]$ modulo $s$ of $Q^*(\sqrt{n})$.  \cite{R7}\\
Let $E_s$ denote the set consisting of classes of the form
$[a,b,c]~(mod~s)$, $n$ modulo $s$ whereas if $n\equiv i(mod~s)$
for some fixed $i\in \{0,1,..., s-1\}$ and the set consisting of
elements of the form $[a,b,c]$ with $n\equiv i(mod~s)$ is denoted
by $E^i_p$ (or $E^n_s$). Obviously $\cup^{s-1}_{i=1}E^i_s=E_s\quad
\textmd{and} \quad E^i_s\cap E^j_s=\phi~~\textmd{for}~~i\neq j$. \cite{R12}\\
The  classification of the real quadratic irrational numbers by
taking prime modulus is very helpful in studying the modular group
action on the real quadratic fields. Thus it becomes interesting
to determine the proper $G$-subsets of $Q^*(\sqrt{n})$ by taking
the action of $G$ on the set $Q^*(\sqrt{n})$ and hence to find the
$G$-orbits of $Q^*(\sqrt{n})$ for each non square $n$.\\

\section{ Modular group $G$  acting on $Q^*(\sqrt{n})$. }
In \cite{R7}, it was shown that the action of the group on
$Q^*(\sqrt{2})$ is transitive, whereas the action of $G$ on
$Q^*(\sqrt{n})$, $n\neq2$ is intransitive. Specifically, it was
proved with the help of classes $[a,b,c](mod~2^2)$ of the elements
of $Q^*(\sqrt{n})$ that $Q^*(\sqrt{n})$,
$n~{\not\equiv}~2(mod~4)$, has two proper $G$-subsets.\\
 Q. Mushtaq \cite{R9}, In the case of $PSL(2,13)$, showed one $G$-orbit
 of length $13$ in the coset diagram for the natural action of
  $PSL(2,Z)$  on any subset of the real projective
 line. In \cite{R12} it was proved that there exist two proper
$G$-subsets of $Q^*(\sqrt{n})$ when $n\equiv0~(mod~p)$ and four
$G$-subsets of $Q^*(\sqrt{n})$ when $n\equiv0~(mod~pq)$. In the
present studies, with the help of the idea of quadratic residues,
we generalize this result and prove some crucial results which
provide us proper $G$-subsets and $G$-orbits  of $Q^*(\sqrt{n})$. \\
We extend this idea to determine four proper $G$-subsets of
$Q^*(\sqrt{n})$ with $n\equiv0~(mod~2pq)$.\\

\textbf{Lemma 2.1}\\
 Let $n\equiv 0(mod~2pq)$ where $p$ and $q$
are two distinct odd primes, Let $\alpha=\frac{a+\sqrt{n}}{c}\in
Q^*(\sqrt{n})$
  then the sets\\
 $S_1=\{\alpha\in Q^*(\sqrt{n}):(c/pq)=1 ~\textmd{or}~(b/pq)=1~\},$\\
 $S_2=\{\alpha\in Q^*(\sqrt{n}):(c/p)=-1~\textmd{or}~(b/p)=-1~\textmd{with}~(c/q)=1 ~\textmd{or}~(b/q)=1~\},$\\
 $S_3=\{\alpha\in Q^*(\sqrt{n}):(c/p)=1~\textmd{or}~(b/p)=1~\textmd{with}~~(c/q)=-1 ~\textmd{or}~(b/q)=-1~\},$\\
 and $S_4=\{\alpha\in Q^*(\sqrt{n}):(c/p)=-1~\textmd{or}~(b/p)=-1~\textmd{with}~(c/q)=-1 ~\textmd{or}~(b/q)=-1~\}$
are four proper $G$-subsets of $Q^*(\sqrt{n})$.\\

\textbf{Proof.}\\
 Let $\frac{a+\sqrt{n}}{c}\in Q^*(\sqrt{n})$ and
$n\equiv 0(mod~2pq)$, then $a^2-n=bc$ forces that
\begin{equation}\label{z1}
a^2\equiv bc(mod~2pq)
\end{equation}
where $a,b,c$ are belonging to the complete residue system $\{0,1,2,...,\bar{2pq}-1\}$.\\
The congruence (\ref{z1}) implies $a^2\equiv bc(mod~2)$,
$a^2\equiv bc(mod~p)$ and  $a^2\equiv bc(mod~q)$ . Since $1$ is
the only quadratic residue of $2$ and there is no quadratic
non-residue of $2$. Thus by Lemma 1.1  the quadratic residues and
quadratic non residues of $pq$ and $2pq$ are the same. We know
that, if $(t,m)=1$ and $m=2pq$, then the congruence $x^2\equiv
t~(mod~m)$ is solvable and has four incongruent solutions if and
only if $t$ is quadratic residue of $m$  \cite{R1}, and in this
case congruence (\ref{z1}) is solvable and has exactly four
  incongruent solutions. \\
   If $a=b=c$ and each of $a,b,c$ are quadratic residue of $pq$
   then there exist four distinct classes
   $$[a,b,c], [-a,b,c], [a,-b,-c], [-a,-b,-c](mod~pq) $$
Thus for each member $[a,b,c](mod~pq)$, we have four
cases.\\
 Case$(i)$  The classes $[a,b,c](mod~pq)$ with  $(bc/pq)=1$,
 Then all these classes are contained in $S_1$.\\
case$(ii)$  The classes $[a,b,c](mod~pq)$ with $(b/p)=-1$,
$(b/q)=1$,  Then all these classes are contained in
$S_2$.\\
case$(iii)$  The classes $[a,b,c](mod~pq)$ with $(b/p)=1$ ,
$(b/q)=-1$,  Then all these classes are contained in
$S_3$.\\
 case$(iv)$   The classes $[a,b,c](mod~pq)$ with $(b/p)=-1$ , $(b/q)=-1$,
 Then all these classes are contained in $S_4$.\\
As  $x(\alpha)=\frac{-a+\sqrt{n}}{b}=\frac{a_1+\sqrt{n}}{c_1}$,
where $a_1=-a,~~b_1=c,~~c_1=b$ and
$y(\frac{a+\sqrt{n}}{c})=\frac{-a+b+\sqrt{n}}{b}=\frac{a_2+\sqrt{n}}{c_2}$,
where $$a_2=-a+b,~~~b_2=-2a+b+c,~\textmd{and}~c_2=b$$ then by the
congruence (\ref{z1}) we have
\begin{equation}\label{z2}
(-a+b)^2\equiv(-2a+b+c)b(mod~p)
\end{equation}
Since the modular group $PSL(2,Z)$ has the representation
$G=\langle x,y:x^2=y^3=1\rangle$ and every element of $G$ is a
word in the generators $x,y$ of G, to prove that $S_1$ is
invariant under the action of G, it is enough to show that every
element of $S_1$ is mapped onto an element of $S_1$ under $x$ and
$y$. Thus clearly by the congruences (1) and (2) the sets $S_1$,
$S_2$, $S_3$ and $S_4$
 are $G$-subsets of $Q^*(\sqrt{n})$. \quad\quad$\Box$\\

\textbf{Remark 2.2}\\
 Since the quadratic residues and quadratic non residues of $pq$ and $2pq$ are the
same. Therefore the number of $G$-subsets of $Q^*(\sqrt{n})$ when
$n\equiv 0(mod~pq)$ or $n\equiv 0(mod~2pq)$
are same.\\

\textbf{Illustration 2.3}\\
In the coset diagram for $Q^*(\sqrt{15})$ there are four
$G$-orbits namely
$$(\sqrt{15})^G,(-\sqrt{15})^G, (\frac{\sqrt{15}}{3})^G ~~\textmd{and}~~
(\frac{\sqrt{15}}{-3})^G$$ and similarly there are four $G$-orbits
for $Q^*(\sqrt{30})$ namely
$$(\sqrt{30})^G,(-\sqrt{30})^G, (\frac{\sqrt{30}}{2})^G  ~~\textmd{and}~~
(\frac{\sqrt{30}}{-2})^G$$. In the closed path lying in the orbit
$(\sqrt{15})^G$, the transformation
$$(yx)^3(y^2x)(yx)^3$$
fixes $k=\sqrt{15}$, that is $g_1(k)=((yx)^3(y^2x)(yx)^3)(k)= k$. \\
 Let $k$ is an ambiguous number then $x(k)$ is also ambiguous but
 one of the number $y(k)$ or $y^2(k)$ is ambiguous.
 The orientation of edges in the coset diagram is
  associated with the involution $x$ and the small triangles with $y$
  which has order 3.  One of $k$ and $x(k)$ is positive and other is negative but one of
  $k$, $y(k)$ or $y^2(k)$ is negative but other two are positive. We use an
  arrow head on an edge to indicate its direction from negative to
  a positive vertex. The following table shows the details of the
  orbits $\alpha^G$, transformations  which fixes $\alpha$, and
   the ambiguous lengths of each orbit.\\

  \begin{table}[htp1]
\caption{ The Orbits of $\alpha \in Q^*(\sqrt{n})$.}
  \begin{center}
\begin{tabular}{|c|c|c|}
\hline
&& \\  $G$-orbits & Transformations & Ambiguous Length \\
 \hline
 && \\ $(\sqrt{15})^G$ & $(yx)^3(y^2x)(yx)^3$ & 14  \\
 \hline
 &&  \\ $(-\sqrt{15})^G$ & $(yx)^3(y^2x)(yx)^3$ & 14  \\
 \hline
  && \\ $(\frac{\sqrt{15}}{3})^G$ & $(yx)(y^2x)^3(yx)$ & 10  \\
 \hline
  && \\  $(\frac{\sqrt{15}}{-3})^G$ & $(yx)(y^2x)^3(yx)$ & 10  \\
 \hline
 && \\ $(\sqrt{30})^G$ & $(yx)^5(y^2x)^2(yx)^5$ & 24  \\
 \hline
 &&  \\ $(-\sqrt{30})^G$ & $(yx)^5(y^2x)^2(yx)^5$ & 24  \\
 \hline
  && \\ $(\frac{\sqrt{30}}{2})^G$ & $(yx)^2(y^2x)(yx)^2(y^2x)(yx)^2$ & 16  \\
 \hline
  && \\  $(\frac{\sqrt{30}}{-2})^G$ & $(yx)^2(y^2x)(yx)^2(y^2x)(yx)^2$ & 16  \\
 \hline
\end{tabular}
\end{center}
\end{table}

Now we extend this idea when $n\equiv 0(mod~p_1 p_2...p_r)$.\\

 \textbf{Theorem 2.4}\\
Let $n\equiv 0(mod~p_1 p_2...p_r)$, where $p_1, p_2,...p_r$ are
distinct odd primes, then there are exactly $2^r$,  $G$-subsets of
$Q^*(\sqrt{n})$.\\

\textbf{Proof.}\\
 Let $\frac{a+\sqrt{n}}{c}\in Q^*(\sqrt{n})$ and
$n\equiv 0(mod~p_1 p_2...p_r)$ where $p_1, p_2,...p_r$ are
distinct odd primes, then $a^2-n=bc$ gives
\begin{equation}\label{z2}
a^2\equiv bc(mod~p_1 p_2...p_r)
\end{equation}
The congruence (\ref{z2}) implies $a^2\equiv bc(mod~p_1)$,
$$a^2\equiv bc(mod~p_2),...,~\textmd{and}~a^2\equiv bc(mod~p_r)$$.\\
We know that, if $(t,m)=1$ and $m =p_1 p_2...p_r$, then the
congruence $x^2\equiv t~(mod~m)$ is solvable if and only if $t$ is
quadratic residue of $m$ \cite{R1}, and in this case congruence
(\ref{z2}) is solvable and has exactly $2^r$ incongruent
solutions. Since all values of $b$ or $c$ which are quadratic
residues and quadratic non-residues of  $m$  lie in the distinct
$G$-subsets and $m$ is the product of $r$ distinct primes,
Thus consequently we obtain $2^r$, $G$-subsets of $Q^*(\sqrt{n})$. \quad\quad$\Box$\\

 \textbf{Corollary 2.5}\\
Let $n\equiv 0(mod~2p_1 p_2...p_r)$ where $p_1, p_2,...p_r$ are
distinct odd primes, then there are exactly $2^r$,  $G$-subsets of
$Q^*(\sqrt{n})$.\\
\textbf{Proof.}\\
 Let $\frac{a+\sqrt{n}}{c}\in Q^*(\sqrt{n})$ and
$n\equiv 0(mod~2p_1 p_2...p_r)$ where $p_1, p_2,...p_r$ are
distinct odd primes, then $a^2-n=bc$ gives
\begin{equation}\label{z3}
a^2\equiv bc(mod~2p_1 p_2...p_r)
\end{equation}
The congruence (\ref{z3}) implies $a^2\equiv bc(mod~2)$,
$a^2\equiv bc(mod~p_1)$,
$$a^2\equiv bc(mod~p_2),...,~\textmd{and}~a^2\equiv bc(mod~p_r)$$
Since $1$ is the only quadratic residue of $2$ and there is no
quadratic non-residue of $2$. Thus by Lemma 1.1 the quadratic
residues and quadratic non residues of $p_1 p_2...p_r$ and $2p_1
p_2...p_r$ are same. Hence the result follows by the Theorem 2.4. \quad\quad$\Box$\\
\textbf{Remark 2.6}\\
The number of $G$-orbits of $Q^*(\sqrt{n})$ when
$n\equiv 0(mod~p_1 p_2...p_r)$ or $n\equiv 0(mod~2p_1 p_2...p_r)$ are same.\\

\textbf{Illustration 2.7}\\
Take $n=3.5.7.11=1155$, Under the action of $G$ on
$Q^*(\sqrt{1155})$ there are
  sixteen $G$-orbits in the coset diagram for $Q^*(\sqrt{1155})$ namely \\
$$(\sqrt{1155})^G,~~~(\frac{\sqrt{1155}}{-1})^G, ~~~(\frac{\sqrt{1155}}{3})^G,~~~ (\frac{\sqrt{1155}}{-3})^G,$$
$$(\frac{\sqrt{1155}}{5})^G,~~~(\frac{\sqrt{1155}}{-5})^G,~~~
(\frac{\sqrt{1155}}{7})^G,~~~ (\frac{\sqrt{1155}}{-7})^G,$$
$$(\frac{\sqrt{1155}}{11})^G,~~~(\frac{\sqrt{1155}}{-11})^G,~~~
(\frac{\sqrt{1155}}{15})^G,~~~ (\frac{\sqrt{1155}}{-15})^G,$$
$$(\frac{\sqrt{1155}}{21})^G,~~~(\frac{\sqrt{1155}}{-21})^G,~~~
(\frac{\sqrt{1155}}{33})^G,~~~ (\frac{\sqrt{1155}}{-33})^G.$$
  Similarly for  $n=2.3.5.7.11=2310$, Under the action of $G$ on
  $Q^*(\sqrt{2310})$ there are sixteen $G$-orbits in the coset diagram
  for $Q^*(\sqrt{2310})$ namely \\
$$(\sqrt{2310})^G,~~~(\frac{\sqrt{2310}}{-1})^G,~~~ (\frac{\sqrt{2310}}{2})^G,~~~ (\frac{\sqrt{2310}}{-2})^G,$$
$$(\frac{\sqrt{2310}}{5})^G,~~~(\frac{\sqrt{2310}}{-5})^G,~~~
(\frac{\sqrt{2310}}{7})^G, ~~~(\frac{\sqrt{2310}}{-7})^G,$$
$$(\frac{\sqrt{2310}}{10})^G,~~~(\frac{\sqrt{2310}}{-10})^G,~~~
(\frac{\sqrt{2310}}{11})^G, ~~~(\frac{\sqrt{2310}}{-11})^G,$$
$$(\frac{\sqrt{2310}}{14})^G,~~~(\frac{\sqrt{2310}}{-14})^G,~~~
(\frac{\sqrt{2310}}{22})^G, ~~~(\frac{\sqrt{2310}}{-22})^G.$$

\textbf{Theorem 2.8}\\
 Let $h=2k+1\geq3$ then there are exactly
two $G$-orbits of $Q^*(\sqrt{2^h})$ namely $ (2^k\sqrt{2})^G$ and
$(\frac{2^k\sqrt{2}}{-1})^G$.\\
\textbf{Proof.}\\
Let $\frac{a+\sqrt{n}}{c}\in Q^*(\sqrt{2^h})$, then $a^2-n=bc$
forces that\\
$$a^2-2^h\equiv bc(mod~2^h)~~\Rightarrow~~a^2\equiv bc(mod~2^h)$$
But the congruence $a^2\equiv bc(mod~2^h)$ is solvable if and only
if $bc\equiv 1(mod~8)$, Moreover the quadratic residue of $2^h,
h\geq3$ are those integers of the form $8l+1$ which are less than
$2^h$. Since all values of $b$ or $c$ which are quadratic residues
 and quadratic non-residues of $2^h$ lie in the distinct orbits. Thus the classes
$[a,b,c]$ (modulo $2^h$) with $b$ or $c$ quadratic residues of
$2^h$ lie in the orbit $(2^k\sqrt{2})^G$ and similarly
 the classes $[a,b,c]$ (modulo $2^h$) with $b$ or $c$ quadratic
non-residues of $2^h$ lie in the orbit
$(\frac{2^k\sqrt{2}}{-1})^G$, This proves the result. \quad\quad$\Box$\\

\textbf{Illustration 2.9}\\
There are exactly two $G$-orbits of $Q^*(\sqrt{2^7})$  namely
$(2^3\sqrt{2})^G$ and $(\frac{2^3\sqrt{2}}{-1})^G$, In the closed
path lying in the orbit $(2^3\sqrt{2})^G$, the transformation
$$(yx)^{11}(y^2x)^3(yx)^5(y^2x)^{3}(yx)^{11}$$  fixes $2^3\sqrt{2}$.
Similarly in the closed path lying in the orbit
$(\frac{2^3\sqrt{2}}{-1})^G$, the transformation
$$(yx)^{11}(y^2x)^3(yx)^5(y^2x)^{3}(yx)^{11}$$  fixes
$-2^3\sqrt{2}$.\\

\section{ Action of the subgroup $~G^*=~\langle yx \rangle$ and $G^{**}=~\langle yx,y^2x \rangle$ on
$Q^*(\sqrt{n})$. } Let us suppose that  $G^*=\langle yx\rangle$
and $G^{**}=\langle yx,y^2x \rangle$ are two subgroups of $G$. In
this section, we determine the $G$-subsets and $G$-orbits  of
$Q^*(\sqrt{n})$ by subgroup $G^*$ and $G^{**}$ acting on
$Q^*(\sqrt{n})$.  Let $yx(\alpha)=\alpha+1$ and
$y^2x(\alpha)=\frac{\alpha}{\alpha+1}$. Thus
$yx(\frac{a+\sqrt{n}}{c})= \frac{a_1+\sqrt{n}}{c_1}$, with
$a_1=a+c,~~b_1=2a+b+c,~~c_1=c$ and $y^2x(\frac{a+\sqrt{n}}{c})=
\frac{a_2+\sqrt{n}}{c_2}$, with
$a_2=a+b,~~~b_2=b,~\textmd{and}~c_2=2a+b+c$.\\
 In the next Lemma
we see that the transformation $yx$ fixes the classes $[0,0,c]$
(modulo $p$) and the chain of these classes help us in finding
$G^*$-subsets of $Q^*(\sqrt{n})$.\\

\textbf{Lemma 3.1}\\
Let $p$ be any prime, $n\equiv 0(mod~p)$, Then for any $k\geq1$,
$(yx)^k [0,0,c]=[kc,k^2c,c]~(mod~p)$ and in particular $(yx)^p
[0,0,c]=[0,0,c]~(mod~p)$.\\
 \textbf{Proof.}\\
 Let $\alpha= [0,0,c]~(mod~p)$ be a class contained in $E^0_p$.
 Applying the linear fractional transformation $yx$ on
$\alpha$ successively we see $yx[0,0,c]=[c,c,c]$, $(yx)^2
[0,0,c]=[2c,4c,c]$, $(yx)^3 [0,0,c]=[3c,9c,c]$  continuing this
process $k$-times we obtain \\
$(yx)^k [0,0,c]=[kc,k^2c,c]~(mod~p)$\\
 In particular for $k=p$,
$kp\equiv 0(mod~p)$ and $k^2p\equiv
0(mod~p)$. Thus we get $(yx)^p [0,0,c]=[0,0,c]~(mod~p)$. \quad\quad$\Box$\\

\textbf{Lemma 3.2}\\
Let $p$ be an odd prime, $n\equiv 0(mod~p)$ and $G^*=\langle
yx \rangle$, Then the sets\\
 $A_1=\{\alpha\in Q^*(\sqrt{n}):(c/p)=1\}$,
$A_2=\{\alpha\in Q^*(\sqrt{n}):(c/p)=-1\}$,\\
$C_1=\{\alpha\in Q^*(\sqrt{n}):c\equiv
0(mod~p)~\textmd{with}~(b/p)=1\}$, $C_2=\{\alpha\in
Q^*(\sqrt{n}):c\equiv 0(mod~p)~\textmd{with}~(b/p)=-1\}$,\\
are $G^*$-subsets of $Q^*(\sqrt{n})$.\\
\textbf{Proof.}\\
 For any  $\alpha=\frac{a+\sqrt{n}}{c}\in
A_1$ with $n\equiv 0(mod~p)$, $a^2-n= bc$ gives
\begin{equation}\label{z4}
a^2\equiv bc(mod~p)
\end{equation}
we have two cases.\\
(i) If $a\equiv0(mod~p)$, The congruence (\ref{z4}) forces that
 $bc\equiv0(mod~p)$, Then either
$b\equiv0(mod~p)$ or $c\equiv0(mod~p)$ but not both. So in this
case $\alpha$ belongs to
the class $[0,b,0]$ or $[0,0,c]$ modulo $p$.\\
(ii) If $a~{\not\equiv}~0(mod~p)$,  $a^2\equiv bc(mod~p)$, Then
(\ref{z4})  forces that  either both  $b,c$ are quadratic residues
of $p$  or  both quadratic non-residues of $p$.\\
As $yx:[a,b,c]\rightarrow[a+c,2a+b+c,c]$, Then it is clear that
the set $A_1$ is invariant under the action of the mapping $yx$,
So the set $A_1$ is a $G^*$-subset $Q^*(\sqrt{n})$. Similarly the
set $A_2$ is
$G^*$-subset of $Q^*(\sqrt{n})$.\\
Again for any $\alpha=\frac{a+\sqrt{n}}{c}\in C_1$  by congruence
(\ref{z4}) $c\equiv 0(mod~p)\Rightarrow a\equiv 0(mod~p)$, with
$b~{\not\equiv}~0(mod~p)$, so the classes belonging to the set
$C_1$ are of the form $[0,b,0]$ with $b$ quadratic residue of $p$.
Since the mapping $yx$ fixes the classes $[0,b,0]$. Thus clearly
the set $C_1$ is a $G^*$-subsets. Similarly the set $C_2$ is
$G^*$-subsets of
$Q^*(\sqrt{n})$. \quad\quad$\Box$\\
In the next theorem we determine two $G$-subsets of
$Q^*(\sqrt{n})$ by using  $A_1$, $A_2$, $C_1$
and $C_2$ as given in Lemma 3.2.\\

\textbf{Theorem 3.3}\\
 The sets $S_1={A_1\cup C_1}$ and $S_2={A_2\cup C_2}$ are two
 $G$-subsets of $Q^*(\sqrt{n})$.\\
\textbf{Proof.}\\
  Let  $\alpha=\frac{a+\sqrt{n}}{c}\in S_1$ then either
  $\alpha\in A_1$ or $\alpha\in C_1$ with  $n\equiv 0(mod~p)$.
  Thus it is clear that the classes $[a,b,c]~(mod~p)$ with $b$ or $c$ quadratic
  residues of $p$ is contained in ${A_1\cup C_1}$. By Lemma 3.1 $yx$ fixes the classes $[0,0,c]$
(modulo $p$). Also the classes belonging to $A_1$, $A_2$ are
connected to the classes belonging to $C_1$, $C_2$, respectively
under $x$. Since the modular group $PSL(2,Z)$ has the
representation $G=\langle x,y:x^2=y^3=1\rangle$ and every element
of $G$ is a word in its generators $x,y$, to prove that $S_1$ is
invariant under the action of G, it is enough to show that every
element of $S_1$ is mapped onto an element of $S_1$ under $x$ and
$y$. Thus clearly  we see that $S_1={A_1\cup C_1}$
   and $S_2={A_2\cup C_2}$ are both $G$-subsets of
   $Q^*(\sqrt{n})$. \quad\quad$\Box$\\

  In view of the above theorem we observe that for $n=2$ the action of $G$ on
  $Q^*(\sqrt{n})$ is transitive. Since $1$ is the only quadratic residue of $2$ and there is no
quadratic non-residue of $2$, Therefore the set $S_2$ becomes
empty and $S_1$ is the only $G$-subset of $Q^*(\sqrt{2})$. While
the action of $G$ on   $Q^*(\sqrt{n})$, $n\neq 2$ is
intransitive.\\

\textbf{Illustration 3.4}\\
Let $p=5$ and $\alpha=\frac{a+\sqrt{n}}{c}\in Q^*(\sqrt{n})$, with
$n\equiv 0(mod~5)$, is of the form $[a,b,c](mod~5)$. \\
 In modulo  $5$, the squares of the integers $1,2,3,4$ are
$$1^2\equiv 4^2\equiv 1~~ \textmd{and}~~ 2^2\equiv 3^2\equiv 4$$  Consequently, the
quadratic residues of $5$ are $1,4$, and the non residues are
$2,3$. Thus  $A_1$ consists of elements of $Q^*(\sqrt{n})$  of the
form
\begin{eqnarray*}
&&[0,0,1],[0,0,4],[1,1,1],[4,1,1],[2,4,1],[2,1,4],[3,1,4],[3,4,1],[1,4,4],\nonumber
\\&&[4,4,4] \textmd{mod}~5~\textmd.
\end{eqnarray*}
Then $A_1$ is invariant under $yx$, Thus $A_1$ is $G^*$-subset of $Q^*(\sqrt{n})$.\\
The elements of $A_2$ are of the form
\begin{eqnarray*}
&&[0,0,2],[0,0,3],[2,2,2],[3,2,2], [2,3,3],[4,3,2],\nonumber
\\&&[4,2,3],[1,2,3),[1,3,2]~\textmd{and}~[3,3,3]~\textmd{mod}~5~\textmd{only}.
\end{eqnarray*}
Again $A_2$ is invariant under $yx$, Thus $A_2$ is also
$G^*$-subset of $Q^*(\sqrt{n})$. \\
The elements of $C_1$ are of the form
$[0,1,0]~\textmd{and}~[0,4,0]$, and the elements of $C_2$ are of
the form $[0,2,0]~\textmd{and}~[0,3,0]$
Thus  $C_1$  and $C_2$ are $G^*$-subsets.\\
Then clearly the sets $S_1={A_1\cup C_1}$ and $S_2={A_2\cup C_2}$
are two $G$-subsets of $Q^*(\sqrt{n})$.\\

In the next lemma we find the conditions when $n$ is quadratic
residue of $p$.\\

\textbf{Lemma 3.5}\\
For any  $\alpha(a,b,c)\in Q^*(\sqrt{n})$,  $n$ is quadratic
residue of $p$ if and only if either $c~\textmd{or}~
b~{\equiv}~0(mod~p)$.\\
\textbf{Proof.}\\
 Let  $\alpha=\frac{a+\sqrt{n}}{c}\in
Q^*(\sqrt{n})$,  Then $a^2-n= bc$ forces that
\begin{equation}\label{z5}
a^2-n \equiv bc(mod~p)
\end{equation}
Let either $b~{\equiv}~0(mod~p)$ or $c~{\equiv}~0(mod~p)$\\
then congruence (\ref{z5}) implies that $a^2\equiv n(mod~p)$
 which shows that  $n$ is quadratic residue of $p$\\
Conversely let $n$ be quadratic residue of $p$ then clearly
$a^2\equiv n(mod~p)$ that shows either
$b~{\equiv}~0(mod~p)$ or $c~{\equiv}~0(mod~p)$.  \quad\quad$\Box$\\
Further we see the action of $G^{**}=\langle yx,y^2x \rangle$ on
$Q^*(\sqrt{n})$ with $n$ quadratic residue of $p$ and determine
four proper $G^{**}$-subsets of $Q^*(\sqrt{n})$.\\

 \textbf{Theorem 3.6}\\
Let $p$ be an odd prime and $n$ is quadratic residue of $p$, let
 $\alpha=\frac{a+\sqrt{n}}{c}\in
Q^*(\sqrt{n})$ and $G^{**}=\langle
yx,y^2x\rangle$, then the sets \\
 $G_1 =\{\alpha\in
Q^*(\sqrt{n}):(c/p)=1\}$, $G_2 =\{\alpha\in Q^*(\sqrt{n}):(c/p)=-1\}$,\\
are two $G^{**}$-subsets of $Q^*(\sqrt{n})$. \quad\quad$\Box$\\
\textbf{Proof.}\\
Let $\alpha=\frac{a+\sqrt{n}}{c}\in Q^*(\sqrt{n})$, with
$a^2-n\equiv bc(mod~p)$, and $a,b,c$ modulo $p$ are belonging to
the set $\{0,1,2,...,p-1\}$.\\
 Let $p$ is an odd prime with $n$ quadratic residue of $p$, then either
$b~{\equiv}~0(mod~p)$ or $c~{\equiv}~0(mod~p)$.\\
Since  $yx:[a,b,c]\rightarrow[a+c,2a+b+c,c]$ and
$y^2x:[a,b,c]\rightarrow[a+b,b,2a+b+c]$, Since every element of
$G^{**}$ is a word in its generators $yx,y^2x$, Then clearly the
sets  $G_1$,  $G_2$,
are two $G^{**}$-subsets of $Q^*(\sqrt{n})$.  \quad\quad$\Box$\\

\textbf{Corollary 3.7} $\langle G^{**}, x \rangle = G$\\
\textbf{Proof.}\\
We know that  $G=\langle x,y:x^2=y^3=1\rangle$ and $G^{**}=\langle
yx,y^2x \rangle$,  the result follows from the fact that the
generators $x,y$ of G  can be written as word  by the elements
of $\langle G^{**}, x \rangle$ . \quad\quad$\Box$\\
Finally we find $G$-orbits of $Q^*(\sqrt{37})$  with help of
$G^{**}$-subsets as given in Theorem 3.6. It is important to note
that $37$ is the smallest prime which have four $G$-orbits and all
odd primes less than $37$ has exactly
two $G$-orbits. \\

 \textbf{Illustration 3.9}\\
 In the coset diagram for $Q^*(\sqrt{37})$, There
are exactly four $G$-orbits of  $Q^*(\sqrt{37})$ given by \\
 $(\sqrt{37})^G$, $(\frac{1+\sqrt{37}}{2})^G$,
$(\frac{1+\sqrt{37}}{-3})^G$, $(\frac{-1+\sqrt{37}}{-3})^G$.\\
The $G_1$ contains three orbits  $(\sqrt{37})^G$,
 $(\frac{1+\sqrt{37}}{-3})^G$,
$(\frac{-1+\sqrt{37}}{-3})^G$, While the set $G_2$ contains only
one orbit $(\frac{1+\sqrt{37}}{2})^G$.\\
 In the closed path lying in the orbit $(\sqrt{37})^G$, the
transformation
$$g_1=(yx)^6(y^2x)^{12}(yx)^6$$
fixes $k=\sqrt{37}$, that is $g_1(k)=((yx)^6(y^2x)^{12}(yx)^6)(k)=
k$, and so gives the quadratic equation $k^2+37=0$, the zeros,
$\pm \sqrt{37}$, of this equation are fixed points of the
transformations $g_1$. \\
In the closed path lying in the orbit $(\frac{1+\sqrt{37}}{2})^G$,
the transformations
 $$g_2=(yx)^3(y^2x)(yx)(y^2x)^5(yx)(y^2x)(yx)^2$$
  fixes $l=\frac{1+\sqrt{37}}{2}$ and so gives the quadratic equation $l^2-l-9=0$, the
zeros, $\frac{1 \pm\sqrt{37}}{2}$, of this equation are fixed
points of $g_2$. \\
Similarly in the closed path lying in the orbit
$(\frac{1+\sqrt{37}}{-3})^G$ the transformation
$$g_3=(yx)^2(y^2x)^2(yx)(y^2x)^3(yx)^2(y^2x)(yx)$$ fixes $(\frac{1+\sqrt{37}}{-3})$
 and corresponding to the closed lying in the orbit $(\frac{-1+\sqrt{37}}{-3})^G$, the transformation
 $$g_4=(yx)(y^2x)(yx)^2(y^2x)^3(yx)(y^2x)^2(yx)^2$$
 fixes $(\frac{-1+\sqrt{37}}{-3})$. \\
By \cite{R11}, \cite{R6} we see that  $\tau^*(37)=124$,   That is
there are $124$ ambiguous numbers in the coset diagram for
$Q^*(\sqrt{37})$ while the ambiguous length of  the orbits are
$48$, $28$, $24$ and $24$ respectively. \quad\quad$\Box$

\end{document}